\documentclass[10pt]{article}

\input{xypic.tex} \xyoption{all}
\usepackage{amsthm,amssymb,latexsym}
\usepackage{graphics}

\newcommand{\Z}{\mathbb{Z}} 
\newcommand{\p}[1]{{\mathbb{P}^{#1}}} \newcommand{\pn}{{\mathbb{P}^n}}
\newcommand{\PP}{{\mathbb{P}}}
\newcommand{\ox}{{\cal O}_{X}} \newcommand{\picx}{{\rm Pic}(X)}
\newcommand{\im}{{\rm Im}} \newcommand{\Hom}{{\rm Hom}} \newcommand{\rk}{{\rm rk}}
 \newcommand{\opn}{{\cal O}_{\mathbb{P}^n}}
\newcommand{\oqn}{{\cal O}_{Q_n}}

\newtheorem{theorem}{Theorem}

\newtheorem{proposition}[theorem]{Proposition}
\newtheorem{remark}[theorem]{Remark}
\newtheorem{example}[theorem]{Example}

\begin{document}

\title{On the semistability of instanton sheaves
over certain projective varieties}
\author{Marcos Jardim \\ IMECC - UNICAMP \\
Departamento de Matem\'atica \\ Caixa Postal 6065 \\
13083-970 Campinas-SP, Brazil
\\  \\
Rosa M. Mir\'o-Roig \\
Facultat de Matem\`atiques \\ Departament d'Algebra i Geometria \\
Gran Via de les Corts Catalanes 585 \\ 08007 Barcelona, SPAIN }

\maketitle

\begin{abstract}
We show that instanton bundles of rank $r\le 2n-1$, defined as the
cohomology of certain linear monads, on an $n$-dimensional
projective variety with cyclic Picard group are semistable in the
sense of Mumford-Takemoto. Furthermore, we show that rank $r\le n$
linear bundles with nonzero first Chern class over such varieties
are stable. We also show that these bounds are sharp.
\end{abstract}


\baselineskip18pt


\section{Introduction}
Let $X$ be a nonsingular projective variety of dimension $n$ over
an algebraically closed field $\mathbb{F}$ of characteristic zero,
and let ${\cal L}$ denote a very ample invertible sheaf on $X$; let
${\cal L}^{-1}$ denote its inverse.

Given (finite-dimensional) $\mathbb{F}$-vector spaces $V$, $W$ and
$U$, a {\em linear monad} on $X$ is a complex of sheaves
\begin{equation} \label{monad}
M_{\bullet} ~ : ~ 0 \to V\otimes{\cal L}^{-1} \stackrel{\alpha}{\to} W\otimes\ox
\stackrel{\beta}{\to} U\otimes{\cal L} \to 0
\end{equation}
which is exact on the first and last terms, i.e.
$\alpha\in\Hom(V,W)\otimes H^0({\cal L})$ is injective while
$\beta\in\Hom(W,U)\otimes H^0({\cal L})$ is surjective. The
coherent sheaf $E=\ker\beta/\im\alpha$ is called the cohomology of
the monad $M_{\bullet}$. The set:
$$ S = \{ x\in X ~ | ~ \alpha(x)\in\Hom(V,W) ~ {\rm is~not~injective} \} $$
is a subvariety called the {\em degeneration locus} of the monad $M_{\bullet}$.

A torsion-free sheaf $E$ on $X$ is said to be a {\em linear sheaf}
on $X$ if it can be represented as the cohomology of a linear
monad and it is said to be an {\em instanton sheaf} on $X$ if in
addition it has $c_1(E)=0$.

Linear monads and instanton sheaves have been extensively studied
for the case $X=\pn$ during the past 30 years, see for instance
\cite{J,OSS} and the references therein. Buchdahl has studied
monads over arbitrary blow-ups of $\p2$ \cite{B}. In a recent
preprint, Costa and Mir\'o-Roig have initiated the study of linear
monads and locally-free instanton sheaves over smooth quadric
hypersurfaces $Q_n$ within $\PP^{n+1}$ ($n\ge3$) \cite{CMR}. They
have asked whether every such locally free sheaf of rank $n-1$ is
stable (in the sense of Mumford-Takemoto) \cite[Question
5.1]{CMR}.

The main goal of this paper is to give a partial answer to their question
in a more general context, showing that locally-free instanton sheaves of rank
$r\le 2n-1$ on an $n$-dimensional smooth projective variety with cyclic Picard group
are semistable, while locally-free linear sheaves of rank $r\le n$ and $c_1\ne0$ on
such varieties are stable. Furthermore, we also show that the bounds on the rank are
sharp by providing examples of rank $2n$ instanton sheaves and rank $n+1$ linear sheaves
on $\pn$ which are not semistable.

We conclude the paper by studying the semistability of special sheaves on $Q_n$, as introduced
by Costa and Mir\'o-Roig. Theorem \ref{mainthm2} provides a partial answer to Question 5.2 in
\cite{CMR}, showing that every rank $r\le 2n-1$ locally-free special sheaf $E$ on $Q_n$
with $c_1=0$ is semistable, while every rank $r\le n$ locally-free
special sheaf on $Q_n$ with $c_1\ne 0$ is stable.

\paragraph{Acknowledgment.}
The first author is partially supported by the FAEPEX grants
number 1433/04 and 1652/04 and the CNPq grant number
300991/2004-5. The second author is partially supported by the
grant MTM2004-00666.


\section{Instanton sheaves on cyclic varieties}
Note that if $E$ is the cohomology of a linear monad as in (\ref{monad}), then:
$$ \rk(E)=w-v-u ~~ {\rm and } ~~ c_1(E)=(v-u)\cdot \ell $$
where $w=\dim W$, $v=\dim V$, $u=\dim U$ and $\ell=c_1({\cal L})$. Thus any instanton
sheaf $E$ can be represented as the cohomology of a monad of the following type:
\begin{equation}\label{instmon}
0 \to ({\cal L}^{-1})^{\oplus c} \stackrel{\alpha}{\to} \ox^{\oplus r+2c}
\stackrel{\beta}{\to} {\cal L}^{\oplus c} \to 0
\end{equation}
where $r$ is the rank and $c$ is called the {\rm charge} of $E$. It also
follows that the total Chern class of $E$ is given by, in the case $u=v$:
$$ c(E)=\frac{1}{(1-\ell^2)^c}=(1+\ell^2+\ell^4+\cdots)^c ~~. $$

\begin{remark}\rm
For $X=\pn$, instanton sheaves exist for $r\geq n-1$ and all $c$ \cite{J}. For $X$ being
a smooth quadric hypersurface of dimension $n\ge3$, instanton sheaves exist for $r\geq n-1$
and all $c$ \cite{CMR}. It would be very interesting to obtain existence results for a
wider class of varieties.
\end{remark}

A smooth projective variety $X$ is said to be {\em cyclic} if
$\picx=\Z$. Examples of cyclic varieties are projective spaces,
grassmannians and complete intersection subvarieties of dimension
$n\ge 3$ within $P^N$, $N\ge4$. We can assume without loss of
generality that ${\cal L}\cong\ox(l)$ for some $l\ge1$ and
$\omega_X\cong\ox(\lambda)$ for some integer $\lambda$.

\begin{proposition}\label{vanish}
Let $X$ be a smooth projective cyclic variety of dimension $n$
such that $H^p(\ox(k))=0$ for $1\le p\le n-1$ and $\forall k$.
Let $E$ be the linear sheaf given by the cohomology
of the monad:
\begin{equation}\label{lm}
0 \to \ox(-l)^{\oplus a} \stackrel{\alpha}{\longrightarrow} \ox^{\oplus b}
\stackrel{\beta}{\longrightarrow} \ox(l)^{\oplus c} \to 0 ~~.
\end{equation}
Then, we have:
\begin{enumerate}
\item for $n\ge2$, $H^0(E(k))=H^0(E^*(k))=0$ for all $k\leq-1$,
\item for $n\ge3$, $H^1(E(k))=0$ for all $k\le -l-1$,
\item for $n\ge4$, $H^{i}(E(k))=0$ for all $k$ and $2\le i \le n-2$,
\item for $n\ge3$, $H^{n-1}(E(k))=0$ for all $k\ge \lambda+l+1$,
\item for $n\ge2$, $H^n(E(k))=0$ for all $k\ge \lambda +1$,
\end{enumerate}\end{proposition}

It is not hard to see that complete intersection subvarieties of
dimension $\ge 3$ within $P^n$, $n\ge4$ do satisfy the conditions
of the theorem. Note also that cyclic Fano varieties (i.e.
$\lambda\le-1$) also satisfy $H^p(\ox(k))=0$ for $1\le p\le n-1$.
Indeed, Kodaira Vanishing Theorem tells us that:
$$H^{i}(\ox(k))=0 \text{ for  all } i<n \text{ and } k\leq-1; \text{ and }$$
$$H^{i}(\ox(k)\otimes \omega_X)=0  \text{ for all } i>0 \text{ and } k\ge 1.$$
By Serre's duality $H^{i}(\ox(k)\otimes\omega_X)\cong H^{n-i}(\ox(-k))^*$.
So, we conclude that
$$H^0(\ox(k))=0 \text{ for  all }  k\leq-1,$$
$$H^{i}(\ox(k))=0 \text{ for  all }  k \text{ and } 1\leq i \leq n-1, \text{ and }$$
$$H^n(\ox(k))=0 \text{ for  all }  k\geq \lambda+1.$$

\begin{proof}
Assuming that $E$ is the cohomology of the linear monad (\ref{lm}),
let $K=\ker\beta$; it is a locally-free sheaf of rank $b-c$
fitting into the sequences:
\begin{equation} \label{ker1}
0 \to K(k) \to \ox(k)^{\oplus b}
\stackrel{\beta}{\longrightarrow} \ox(k+l)^{\oplus c} \to 0 ~~
{\rm and}
\end{equation}
\begin{equation} \label{ker2}
0 \to \ox(k-l)^{\oplus a} \stackrel{\alpha}{\longrightarrow}
K(k) \to E(k) \to 0 ~~.
\end{equation}
Passing to cohomology, the exact sequence (\ref{ker1}) yields,
in the appropriate ranges of $n$:
$$ H^0(K(t))=0 \text{ for } t\leq -1~, $$
$$ H^1(K(t))=0 \text{ for } t\le -l-1~, $$
$$ H^{i}(X,K(t))=0 \text{ for all } t \text{ and } 2\le i \le n-1~, $$
$$ H^n(K(t))=0 \text{ for } t\geq \lambda+1~. $$
Passing to cohomology, the exact sequence (\ref{ker2}) yields:
$$ H^0(E(k))=0 \text{ for all } k\leq-1~, $$
$$ H^1(E(k))=0 \text{ for all } k\le -l-1~, $$
$$ H^{i}(E(k))=0 \text{ for all } k \text{ and } 2\le i \le n-2~, $$
$$ H^{n-1}(E(k))=0 \text{ for all } k\ge \lambda +l+1~, $$
$$ H^{n}(E(k))=0 \text{ for all } k\ge \lambda +1~, $$
as desired.

Dualizing  sequences (\ref{ker1}) and (\ref{ker2}), we obtain:
\begin{equation} \label{ker1d}
0 \to \ox(-k-l)^{\oplus c} \stackrel{\beta^*}{\longrightarrow}
\ox(-k)^{\oplus b} \to K^*(-k)  \to 0 ~~ {\rm and}
\end{equation}
\begin{equation} \label{ker2d}
0 \to E^*(-k) \to K^*(-k)
\stackrel{\alpha^*}{\longrightarrow} \ox(-k+l)^{\oplus a} \to
{\cal E}xt^1(E(k),\ox) \to 0 ~~.
\end{equation}
Again, passing to cohomology, (\ref{ker2d}) forces
$H^0(E^*(k))\subseteq H^0(K^*(k))$ for all $k$, while
(\ref{ker1d}) implies $H^0(K^*(k))=0$ for $k\leq-1$. Therefore,
$H^0(E^*(k))=0$ for all $k\le -1$.
\end{proof}

\begin{remark}\label{rem3}\rm
It follows from (\ref{ker2d}) that ${\cal E}xt^1(E,\ox)={\rm coker}\alpha^*$,
i.e. the degeneration locus of the monad (\ref{instmon}) coincides with the
support of ${\cal E}xt^1(E,\ox)$. Furthermore, it also follows from
(\ref{ker2d}) that ${\cal E}xt^p(E,\ox)=0$ for $p\geq 2$.
\end{remark}

\begin{proposition}
Let $E$ be a linear sheaf on a smooth projective variety $X$ (not necessarily cyclic).
\begin{enumerate}
\item $E$ is locally-free if and only if its degeneration locus is empty;
\item $E$ is reflexive if and only if its degeneration locus is a subvariety
of codimension at least 3;
\item $E$ is torsion-free if and only if its degeneration locus is a subvariety
of codimension at least 2.
\end{enumerate} \end{proposition}
\begin{proof}
Let $S$ be the degeneration locus of the linear monad associated
to the linear sheaf $E$. From Remark \ref{rem3}, we know that
${\cal E}xt^p(E,\ox)=0$ for $p\geq2$ and
$$ S = {\rm supp}~{\cal E}xt^1(E,\ox) =
\{ x\in X ~ | ~ \alpha(x) ~ {\rm is~not~injective} ~ \} .$$

The first statement is clear; so it is now enough to argue that
$E$ is torsion-free if and only if $S$ has codimension at least 2
and that $E$ is reflexive if and only if $S$ has codimension at
least 3.

Recall that the $m^{\rm th}$-singularity set of a coherent sheaf $\cal F$
on $X$ is given by:
$$ S_m({\cal F}) = \{ x\in X ~|~ dh({\cal F}_x) \geq n-m \} $$
where $dh({\cal F}_x)$ stands for the homological dimension of
${\cal F}_x$ as an ${\cal O}_x$-module:
$$ dh({\cal F}_x) = d ~~~ \Longleftrightarrow ~~~
\left\{ \begin{array}{l}
{\rm Ext}^d_{{\cal O}_x}({\cal F}_x,{\cal O}_x) \neq 0 \\
{\rm Ext}^p_{{\cal O}_x}({\cal F}_x,{\cal O}_x) = 0 ~~ \forall
p>d. \end{array} \right. $$

In the case at hand, we have that $dh(E_x) = 1$ if $x\in S$, and
$dh(E_x) = 0$ if $x\notin S$. Therefore
$S_0(E)=\cdots=S_{n-2}(E)=\emptyset$, while $S_{n-1}(E)=S$. It
follows that \cite[Proposition 1.20]{ST}:
\begin{itemize}
\item if ${\rm codim}~S\geq 2$, then $\dim S_m(E)\leq m-1$ for all
$m<n$, hence $E$ is a locally 1$^{\rm st}$-syzygy sheaf; \item if
${\rm codim}~S \geq 3$, then $\dim S_m(E)\leq m-2$ for all $m<n$,
hence $E$ is a locally 2$^{\rm nd}$-syzygy sheaf.
\end{itemize}
The desired statements follow from the observation that $E$ is
torsion-free if and only if it is a locally 1$^{\rm st}$-syzygy sheaf,
while $E$ is reflexive if and only if it is a locally
2$^{\rm nd}$-syzygy sheaf \cite[p. 148-149]{OSS}.
\end{proof}

\begin{remark}\rm
Note that if $E$ is a locally-free linear sheaf on $X$, which is
represented as the cohomology of the linear monad
$$ M_{\bullet} ~ : ~ \quad
0 \to ({\cal L}^{-1})^{\oplus a} \stackrel{\alpha}{\to}
\ox^{\oplus b} \stackrel{\beta}{\to} {\cal L}^{\oplus c} \to 0 ~~,$$
its dual $E^*$ is also a linear sheaf, being represented as the cohomology
of the dual monad
$$ M_{\bullet}^* ~ : ~ \quad
0 \to ({\cal L}^{-1})^{\oplus c} \stackrel{\alpha}{\to}
\ox^{\oplus b} \stackrel{\beta}{\to} {\cal L}^{\oplus a} \to 0 ~~.$$
In particular, if $E$ is a locally-free instanton sheaf on $X$ then
its dual $E^*$ is also an instanton. In general, however, there are
non-locally-free instanton sheaves whose duals are not instantons; the
simplest example of this situation is a non-locally-free nullcorrelation
bundle on $\p3$.
\end{remark}


\section{Semistability of instanton sheaves}

Fixed an ample invertible sheaf $\cal L$ with $c_1(\cal L)=\ell$
on a projective variety $X$ of dimension $n$, recall that the
slope $\mu (E)$ with respect to $\cal L$ of a torsion-free sheaf
$E$ on $X$ is defined as follows:
$$ \mu (E):=\frac{c_1(E)\ell ^{n-1}}{rk(E)} ~~. $$

We say that $E$ is semistable with respect to $\cal L$
if for every coherent sheaf $0\ne F\hookrightarrow E$ we have
$\mu (F) \le \mu (E)$. Furthermore, if for every coherent
sheaf $0\ne F\hookrightarrow E$ with $0<rk(F)<rk(E)$ we have
$\mu (F) < \mu (E)$ then $E$ is said to be stable. A sheaf
$E$ is said to be properly semistable if it is semistable but not
stable. It is also important to recall that $E$ is (semi)stable if
and only if $E^*$ is (semi)stable if and only if
$E\otimes {\cal L}^{\otimes k}$ is (semi)stable.

The goal of this section is to study the
(semi)stability of instanton sheaves.

\begin{proposition}\label{rk2}
Every rank 2 torsion-free instanton sheaf on a cyclic variety is
semistable.
\end{proposition}
\begin{proof}
Let us first consider a rank 2 reflexive sheaf $F$ on $X$ with
$c_1(F)=0$ and $H^0(F(-1))=0$; we argue that $F$ is semistable.
Indeed, if $F$ is not semistable, then any destabilizing sheaf
$L\hookrightarrow F$ with torsion-free quotient $F/L$ must be
reflexive (see \cite[p. 158]{OSS}). But every rank 1 reflexive
sheaf is locally-free, thus $L=\ox(d)$ with $d=c_1(L)>0$ since
$\picx=\Z$. It follows that $H^0(F(-d))\neq0$, hence
$H^0(F(-1))\neq 0$ as well.

Now if $E$ is a rank 2 torsion-free sheaf with $c_1(E)=0$ and
$H^0(E^*(-1))=0$, then $F=E^*$ is a rank 2 reflexive sheaf with
$c_1F=0$ and $H^0(F(-1))=0$. But we've seen that such $F$ is
semistable, hence $E$ is also semistable. Together with the first
statement in Proposition \ref{vanish}, the desired result follows.
\end{proof}

For instanton sheaves of higher rank, we have our first main result:

\begin{theorem}\label{mainthm}
Let $E$ be a rank $r$ instanton sheaf on a cyclic variety $X$ of
dimension $n$. If  $E$ is locally-free and $r\le 2n-1$, then $E$
is semistable.
\end{theorem}

Since smooth quadric hypersurfaces are cyclic, the above statement
provides in particular a partial answer to the questions raised in
\cite[Questions 5.1 and 5.2]{CMR}.

The proof of Theorem \ref{mainthm} is based on a very useful to
decide whether a locally-free sheaf on cyclic variety is
(semi)stable.

Recall that for any rank $r$ locally-free sheaf $E$ on a cyclic
variety $X$, there is a uniquely determined integer $k_E$ such
that $-r+1\le c_1(E(k_E))\le 0$. We set $E_{norm}:=E(k_E)$ and we
call $E$ normalized if $E=E_{norm}$. We then have the following
criterion.

\begin{proposition}\label{hoppe}
Let $E$ be a rank $r$ locally-free sheaf on a cyclic variety $X$.
If $H^0((\wedge ^qE)_{\rm norm})=0$ for $1\leq q\leq r-1$, then
$E$ is stable. If $H^0((\wedge ^qE)_{\rm norm}(-1))=0$ for $1\leq
q\leq r-1$, then $E$ is semistable.
\end{proposition}

\begin{proof} \cite{H}; Lemma 2.6.
\end{proof}

\noindent{\em Proof of Theorem \ref{mainthm}.} We argue that every
instanton sheaf on an $n$-dimensional cyclic variety $X$
satisfying the conditions of the theorem fulfills Hoppe's
criterion (see Proposition \ref{hoppe}).

Indeed, let $E$ be a rank $r$ locally-free instanton sheaf on $X$.
Assume that $E$ can be represented as the cohomology of the linear
monad as in (\ref{lm}).

Considering the long exact sequence of symmetric powers associated to the
sheaf sequence
$$ 0 \to K \to \ox^{\oplus r+2c} \to \ox(l)^{\oplus c} \to 0 $$
twisted by $\ox(k)$, we have:
$$ 0 \to \wedge^qK(k) \to \wedge^q(\ox^{\oplus r+2c})(k) \to \cdots $$
$$ \cdots \to \ox^{\oplus r+2c}(k)\otimes S^{q-1}(\ox(l)^{\oplus c})
\to S^q(\ox(l)^{\oplus c}) \to 0 $$
Cutting into short exact sequences, passing to cohomology and using the
fact that $H^p(\ox(k))=0$ for $p\le n-1$ and $k\le-1$ (Kodaira vanishing theorem),
we conclude that
$$ H^p(\wedge^qK(k))=0 ~~{\rm for}~~ 1\leq q \leq r+c={\rm rk}(K)
~~,~~ p\le n-1 ~~{\rm and}~~  k\le -pl-1 .$$

Now consider the long exact sequence of exterior powers
associated to the sheaf sequence
$$ 0 \to \ox(-l)^{\oplus c} \to K \to E \to 0 $$
and twisted by $\ox(-1)$:
$$ 0\to \ox(-ql-1)^{\oplus{c+q-1\choose q}} \to
K((-q+1)l-1)^{\oplus{c+q-2\choose q-1}} \to \cdots $$
\begin{equation}
\cdots \to  \wedge ^{q-1}K(-1-l)^{\oplus c} \to \wedge ^qK(-1) \to
\wedge ^qE(-1) \to 0 ~~.
\end{equation}
Cutting into short exact sequences and passing to cohomology, we
obtain that
\begin{equation}
H^0(\wedge^pE(-1))=0~~ {\rm for} ~~ 1\leq p\leq n-1 ~~.
\end{equation}
If ${\rm rk}(E)\leq n$, this proves that $E$ is semistable by Proposition \ref{hoppe}.

If ${\rm rk}(E)=n+1$, we have, since $c_1(E)=0$ and $E$ is normal:
\begin{equation} \label{h0refl}
H^0(\wedge^nE(-1)) \simeq H^0(E^*(-1))=0 ~~,
\end{equation}
thus $E$ is also semistable.

Assume ${\rm rk}(E)>n+1$. The dual $E^*$ is also a locally-free
instanton sheaf on $X$, so
\begin{equation}\label{dois}
H^0(\wedge^q(E^*)(-1))=0 ~~ {\rm for} ~~ 1\leq q\leq n-1 ~~.
\end{equation}
But $\wedge^p(E)\simeq\wedge^{r-p}(E^*)$, since $\det(E)=\ox$; it
follows that:
\begin{eqnarray}
\nonumber H^0(\wedge^pE(-1))=H^0(\wedge^{r-p}(E^*)(-1))=0 &~~ {\rm
for} ~~&
1\leq r-p \leq n-1 \\
\label{tres}&\Longrightarrow& r-n+1\leq p \leq r-1 ~~.
\end{eqnarray}

Together, (\ref{dois}) and (\ref{tres}) imply that if $E$ is a
rank $r\leq 2n-1$ locally-free instanton sheaf, then:
$$ H^0(\wedge^pE(-1))=0~~ {\rm for} ~~ 1\leq p\leq 2n-2 $$
hence $E$ is semistable by Proposition \ref{hoppe}.
\hfill$\Box$

On the other hand, we have:

\begin{proposition}
Let $H=h^0({\cal L})$. For $r>(H-2)c$, there are no stable rank
$r$ instanton sheaves of charge $c$ on $X$.
\end{proposition}

In particular, for $X=\pn$ and ${\cal L}=\opn(1)$, it follows that
every locally-free instanton sheaf on $\pn$ of charge $1$ and rank
$r$ with $n\le r\le 2n-1$ must be properly semistable; for $X=Q_n$
and ${\cal L}=\oqn(1)$, every locally-free instanton sheaf on
$Q_n$ of charge $1$ and rank $r$ with $n+1\le r\le 2n-1$ must be
properly semistable.

\begin{proof}
For the second part, note that if $E$ is a stable torsion-free sheaf with
$c_1(E)=0$, then $H^0(E)=0$. Indeed, if $H^0(E)\neq0$, then there is a map
$\ox\to E$, which contradicts stability.

It follows from the sequences (\ref{ker1}) and (\ref{ker2}) for $k=0$ that:
$$ H^0(E) \simeq H^0(K) \simeq \ker
\{~ H^0\beta ~:~ H^0(\ox^{\oplus r+2c}) \to H^0({\cal L}^{\oplus c}) ~\} ~~.$$
If $r>(H-2)c$, then
the map $H^0\beta$ cannot be injective, $H^0(E)\neq0$ and $E$ cannot be
stable.
\end{proof}

Now dropping the $c_1(E)=0$ condition, we obtain:

\begin{theorem}\label{mthm2}
Let $E$ be a rank $r\le n$ linear sheaf on a cyclic variety $X$ of dimension $n$.
 If $E$  is locally-free and $c_1(E)\ne 0$, then $E$ is stable.

\end{theorem}

\begin{proof}
Since $E$ is a linear sheaf, it is represented as the cohomology
of a linear monad
$$ 0 \to \ox (-1)^{\oplus a} \stackrel{\alpha}{\to}
\ox^{\oplus b} \stackrel{\beta}{\to} \ox (1)^{\oplus c} \to 0
~~,$$ so that $c_1(E)=(a-c)\ell$.

Assuming $a-c>0$, we have $\mu (\wedge ^qE)=q(a-c)/r>0$, hence
$(\wedge^qE)_{\rm norm}=(\wedge ^qE)(t)$ for some $t\le -1$.

On the other hand, arguing as in the proof of Theorem \ref{mainthm} we get
\begin{equation}\label{firsthalf}
H^0((\wedge^qE)(-1)) = 0 \text{ for all } q\le n-1~~.
\end{equation}

Therefore, if $E$ is a rank $r\leq n$ locally-free sheaf
represented as the cohomology of a linear monad and $c_1(E)> 0$,
then:
$$ H^0((\wedge^pE)_{norm})=0~~ {\rm for} ~~ 1\leq p\leq r-1. $$
Hence $E$ is stable by Proposition \ref{hoppe}.

For the second statement, note that if $E$ is a locally-free linear
sheaf with $c_1(E)<0$, then $E^*$ is a locally-free linear
sheaf with $c_1(E^*)>0$. By the argument above, $E^*$ is stable; hence
$E$ is stable whenever $c_1(E)\ne0$, as desired.
\end{proof}

\vskip 2mm

We will end this section with examples which
illustrate that the upper bounds in the rank given in Theorems
\ref{mainthm} and \ref{mthm2} are sharp.
To establish them, we first need to provide the following useful
cohomological characterization of linear sheaves on projective
spaces.

\begin{proposition}\label{criterio}
Let $F$ be a torsion-free sheaf on $\pn$. $F$ is the cohomology of a
linear monad
$$ 0 \to \opn(-1)^{\oplus a} \stackrel{\alpha}{\longrightarrow} \opn^{\oplus b}
\stackrel{\beta}{\longrightarrow} \opn(1)^{\oplus c} \to 0 $$
if and only if the following cohomological conditions hold:
\begin{itemize}
\item for $n\ge2$, $H^0(F(-1))=0$ and $H^n(F(-n))=0$;
\item for $n\ge3$, $H^1(F(k))=0$ for $k\leq-2$ and $H^{n-1}(F(k))=0$ for $k\geq-n+1$;
\item for $n\ge4$, $H^p(F(k))=0$ for $2\le p\le n-2$ and all $k$.
\end{itemize}
\end{proposition}

\begin{proof}
The fact that linear sheaves satisfy the cohomological conditions above is a
consequence of Proposition \ref{vanish}.

For the converse statement, first note that $H^0(F(-1))=0$ implies that
$H^0(F(k))=0$ for $k\leq-1$, while $H^n(F(-n))=0$ implies that $H^n(F(k))=0$
for $k\ge-n$. Moreover, we claim that ($q=0,\dots,n$ and $p=0,-1,\dots,-n$):
\begin{equation}\label{vhq}
H^q(F(-1)\otimes\Omega_{\pn}^{-p}(-p))=0 ~~ {\rm for}~q\neq1 ~~
{\rm and~for}~ q=1,~p\leq-3 ~ .
\end{equation}

Now the key ingredient is the {\em Beilinson spectral sequence}
\cite{OSS}: for any coherent sheaf $G$ on $\pn$, there exists a
spectral sequence $\{E^{p,q}_r\}$ whose $E_1$-term is given by
($q=0,\dots,n$ and $p=0,-1,\dots,-n$):
$$ E_1^{p,q} =  H^q(G\otimes\Omega_{\pn}^{-p}(-p))\otimes \opn(p) $$
which converges to
$$ E^i = \left\{ \begin{array}{l}
G ~,~ {\rm if} ~ p+q=0 \\ 0 ~ {\rm otherwise} \end{array} \right.
~~ . $$ Applying the Beilinson spectral sequence to $G=F(-1)$, it
then follows that it degenerates at the $E_2$-term, so that the
monad
\begin{eqnarray}
\label{m3} 0 & \to & H^1(F(-1)\otimes\Omega_{\pn}^2(2))\otimes\opn(-2) \to \\
\nonumber & \to & H^1(F(-1)\otimes\Omega^1_{\pn}(1))\otimes\opn(-1) \to
H^1(F(-1))\otimes\opn \to 0
\end{eqnarray}
has $F(-1)$ as its cohomology. Tensoring (\ref{m3}) by $\opn(1)$, we conclude
that $F$ is the cohomology of a linear monad, as desired.

The claim (\ref{vhq}) follows from repeated use of the exact sequence
$$ H^q(F(k))^{\oplus m} \to H^q(F(k+1)\otimes\Omega_{\pn}^{-p-1}(-p-1)) \to $$
\begin{equation} \label{e(k)euler hom}
\to H^{q+1}(F(k)\otimes\Omega_{\pn}^{-p}(-p)) \to H^{q+1}(F(k))^{\oplus m}
\end{equation}
associated with Euler sequence for $p$-forms on $\pn$ twisted by $F(k)$:
\begin{equation} \label{e(k)euler}
0 \to F(k)\otimes\Omega_{\pn}^{-p}(-p) \to F(k)^{\oplus m} \to
F(k)\otimes\Omega_{\pn}^{-p-1}(-p) \to 0 ~,
\end{equation}
where $q=0,\dots,n$ , $p=0,-1,\dots,-n$ and
$m=\left(\begin{array}{c}n+1\\-p\end{array}\right)$.
\end{proof}

We are finally ready to construct rank $2n$ locally-free instanton
sheaves on $\pn$ which are not semistable; in other words the
bound $r\le 2n-1$ in the second part of Theorem \ref{mainthm} is sharp.

\begin{example}\label{ex1} \rm
Let $X=\pn$, $n\ge 4$. By Fl\o ystad's theorem \cite{F}, there is
a linear monad:
\begin{equation} \label{def-F}
0 \to \opn(-1)^{\oplus 2} \stackrel{\alpha}{\to} \opn^{\oplus n+3}
\stackrel{\beta}{\to} \opn(1) \to 0
\end{equation}
whose cohomology $F$ is a locally-free sheaf of rank $n$ on $\pn$
and $c_1(F)=1$.

Dualizing we get a linear monad:
$$ 0 \to \opn(-1) \stackrel{\beta^*}{\to} \opn^{\oplus n+3}
\stackrel{\alpha^*}{\to} \opn(1)^{\oplus 2} \to 0 $$ whose
cohomology is $F^*$, hence it is a locally-free linear sheaf of rank
$n$ on $\pn$ and $c_1(F^*)=-1$.

Take an extension $E$ of $F^*$ by $F$:
$$ 0\to F\to E\to F^*\to 0. $$
Using the cohomological criterion given in Proposition \ref{criterio},
it is easy to see that the extension of linear sheaves is also a linear
sheaf. Moreover, $c_1(E)=0$, i.e. $E$ is a rank $2n$ locally-free
instanton sheaf of charge 3 which is not semistable.

Such extensions are classified by ${\rm Ext}^1(F^*,F)=H^1(F\otimes F)$.
We claim that there are non-trivial extensions of $F^*$ by $F$.
Indeed, we consider the exact sequences
\begin{equation} \label{seq1}
0 \to K=\ker(\beta)\to \opn^{\oplus n+3} \stackrel{\beta}{\to} \opn(1) \to 0 ~,
\end{equation}
\begin{equation} \label{seq2}
0 \to \opn(-1)^{\oplus 2} \to K \to F \to 0
\end{equation}
associated to the linear  monad (\ref{def-F}). We apply the exact covariant functor
$\cdot \otimes F$ to the exact sequences (\ref{seq1}) and (\ref{seq2}) and we
obtain the exact sequences
$$ 0 \to K\otimes F\to F^{\oplus n+3} \to F(1) \to 0 ~, $$
$$ 0 \to F(-1)^{\oplus 2} \to K\otimes F \to F \otimes F \to 0 ~~. $$
Using Proposition \ref{vanish}, we obtain
$H^{i}(K\otimes F)=H^{i}(F\otimes F)=0$ for all $i\ge 3$. Hence,
$\chi (F\otimes F)= h^0((F\otimes F))-h^1((F\otimes F))+h^2((F\otimes F))$.
On the other hand,
$$\chi (F\otimes F)= \chi(K\otimes F)-2\chi(F(-1))=$$
$$(n+3)\chi(F)-\chi(F(1))-2\chi(F(-1))=8-\frac{n^2}{2}-\frac{n}{2}<0 ~~,~~
{\rm if}~~ n\ge4 ~~. $$
Thus if $n\ge4$, we must have $h^1((F\otimes F))>0$, hence
there are non-trivial extensions of $F^*$ by $F$.

For $X=\pn$, $2\le n \le 3$, arguing as above, we can
construct a rank $2n$ locally-free instanton which is not
semistable as a non-trivial extension $E$ of $F^*$ by $F$, where
$F$ is a linear sheaf represented as the cohomology of the linear
monad
$$ 0 \to \opn(-1)^{\oplus 4} \stackrel{\alpha}{\to} \opn^{\oplus n+7}
\stackrel{\beta}{\to} \opn(1)^{\oplus 3} \to 0 .$$
\end{example}

%
%
%

To conclude this section, we show that the upper bound in the rank given in
Theorem \ref{mthm2} is also sharp:

\begin{example}\label{ex2} \rm
Let $X=\pn$, $n\ge2$. By Fl\o ystad's theorem \cite{F}, there is a
linear monad:
\begin{equation}\label{mon.g}
0 \to \opn(-1)^{\oplus 4} \stackrel{\alpha}{\to} \opn^{\oplus n+9}
\stackrel{\beta}{\to} \opn(1)^{\oplus 5} \to 0
\end{equation}
whose cohomology $G$ is a locally-free sheaf of rank $n$ on $\pn$ and
$c_1(G)=-1$.

Now $G^*$ is the cohomology of the dual monad
$$ 0 \to \opn(-1)^{\oplus 5} \stackrel{\beta^*}{\to} \opn^{\oplus n+9}
\stackrel{\alpha^*}{\to} \opn(1)^{\oplus 4} \to 0 ~~. $$
It follows that:
$$ H^1(G^*)=H^1(\ker\alpha^*)=
{\rm coker}\{H^0\alpha^*:H^0(\opn^{\oplus n+9}) \to
H^0(\opn(1)^{\oplus 4})\} ~~.$$ Since $n\ge2$ forces $4n+4>n+9$,
the generic map $\alpha$ will have ${\rm coker}(H^0\alpha^*)\ne0$.
In other words, there exists a rank $n$ locally-free linear sheaf
$G$ on $\pn$ with $c_1(G)=-1$ and $H^1(G^*)\neq0$.

Take an extension $E$ of such a linear sheaf $G$ by $\opn$:
\begin{equation}\label{e2}
0 \to \opn \to E \to G \to 0.
\end{equation}
Using the cohomological criterion given in Proposition
\ref{criterio}, it is easy to see that $E$ is a rank $n+1$
locally-free linear sheaf with $c_1(E)=c_1(G)=-1$. It is not
stable, since $H^0(E)\ne0$.

Note also that there are nontrivial extensions of $G$ by $\opn$
since $H^1(G^*)\neq0$. Furthermore, the dual $E^*$ is an example
of a rank $n+1$ locally-free linear sheaf with $c_1(E)>0$ which is
not stable.
\end{example}

We do not know how to estabilish the semistability of torsion-free
instanton sheaves of rank higher than $3$. However, for each $n\ge2$,
it is easy to show, using the same technique as in the examples above,
that there are unstable torsion-free instanton sheaves of rank $n+1$
in $\pn$, see \cite[Example 3]{J}. The natural, sharp conjecture would be
that every torsion-free instanton sheaf of rank $r\leq n$ on a cyclic
variety $X$ of dimension $n$ is semistable; this statement is true for
$n=2$, see Proposition \ref{rk2} above.

For reflexive linear sheaves, one can construct rank $n+2$ reflexive instanton
sheaves which are not semistable; in this case, one can expect that every
reflexive instanton sheaf of rank $r\le n+1$ on a cyclic variety $X$ of dimension
$n$ is semistable.


\section{Special sheaves on smooth quadric \\ hypersurfaces}

Now we restrict ourselves to the set-up in \cite{CMR}, and we assume
that $Q_n$ is a smooth quadric hypersurface within $\mathbb{P}^{n+1}$,
$n\geq 3$; such varieties are cyclic.

Recall that a {\em special sheaf} $E$ on $Q_n$ is defined
\cite[Definition 3.4]{CMR} as either the cohomology of a linear monad
$$ {\rm (M1)} ~~ 0 \to \oqn(-1)^{\oplus a} \to \oqn^{\oplus b}
\to \oqn(1)^{\oplus c} \to  0 ~~, $$
or the cohomology of a  monad of the following type
$$ {\rm (M2.1)} ~~ 0 \to \Sigma(-1)^{\oplus a} \to \oqn^{\oplus b} \to
\oqn(1)^{\oplus c} \to  0 ~~, \text{ if } n \text{ is odd, } $$
$$ {\rm (M2.2)} ~~ 0 \to \Sigma_1(-1)^{\oplus a_1}\oplus \Sigma_2(-1)^{\oplus a_2}
\to \oqn^{\oplus b} \to \oqn(1)^{\oplus c} \to  0 ~~, \text{ if }
n \text{ is even, } $$ where $\Sigma$ is the Spinor bundle for $n$
odd, and $\Sigma_1,\Sigma_2$ are the Spinor bundles for $n$ even.

Clearly, instanton sheaves on $Q_n$ are special sheaves of the first
kind with zero degree.

\begin{proposition}\label{vanish2}
Let $E$ be a  special sheaf on $Q_n$, $n\ge3$. Then one of the
following conditions holds:
\begin{enumerate}
\item $E$ is the cohomology of a linear monad, and
\begin{itemize}
\item $H^0(E(k))=H^0(E^*(k))=0$ for all $k\leq-1$, \item
$H^1(E(k))=0$ for all $k\le -2$, \item $H^{i}(E(k))=0$ for all $k$
and $2\le i \le n-2$, \item $H^{n-1}(E(k))=0$ for all $k\ge -n
+2$, \item $H^n(E(k))=0$ for all $k\ge -n+1$,

and if $E$ is locally-free: \item $H^n(E^*(k))=0$ for all $k\ge
-n+1$; or
\end{itemize}

\item $E$ is the cohomology of a monad of type
{\em (M2.1)} and {\em (M2.2)}, and
\begin{itemize}
\item $H^0(E(k))=H^0(E^*(k))=0$ for all $k\leq-1$, \item
$H^1(E(k))=0$ for all $k\le -2$, \item $H^{i}(E(k))=0$ for all $k$
and $2\le i \le n-2$, \item $H^{n-1}(E(k))=0$ for all $k\ge -n
+1$, \item $H^n(E(k))=0$ for all $k\ge -n+1$,

 and if $E$ is locally-free: \item $H^n(E^*(k))=0$ for all $k\ge -n+1$
\end{itemize} \end{enumerate}
\end{proposition}

\begin{proof}
(1) It is analogous to the proof of Proposition
\ref{vanish}.

(2) If $n$ is odd we consider the exact sequences $$  0 \to
ker(\delta)  \to \oqn^{\oplus b} \stackrel{\delta}{\to}
\oqn(1)^{\oplus c} \to 0 ~~, $$ $$  0 \to \Sigma(-1)^{\oplus a}
\to ker(\delta)  \to E \to 0
$$  and if $n$ is even we consider the exact sequences $$  0 \to
ker(\psi )  \to \oqn^{\oplus b} \stackrel{\psi }{\to}
\oqn(1)^{\oplus c} \to 0 ~~, $$ $$  0 \to \Sigma_1(-1)^{\oplus
a_1}\oplus \Sigma_1(-2)^{\oplus a_2} \to ker(\psi) \to E \to 0
$$ and we argue as in the proof of Proposition \ref{vanish} taking
into account that

$H^0(\Sigma (k))=H^0(\Sigma _1(k))=H^0(\Sigma _2(k))=0$ for all $k\le
-1$,

$H^{i}(\Sigma (k))=H^{i}(\Sigma _1(k))=H^{i}(\Sigma _2(k))=0$ for all
$k$ and $1\le i \le n-1$, and

$H^{n}(\Sigma (k))=H^{n}(\Sigma _1(k))=H^{n}(\Sigma _2(k))=0$ for all
$k\ge n.$
\end{proof}

\begin{proposition}
Every rank 2 torsion-free special sheaf $E$ on $Q_n$ with $c_1(E)=0$
is semistable.
\end{proposition}
\begin{proof}
Since every torsion-free special sheaf $E$ on $Q_n$ satisfies
$H^0(E(k))=H^0(E^*(k))=0$, simply use the argument in the proof of
Proposition \ref{rk2}.
\end{proof}

Finally, for higher rank locally-free special sheaves on $Q_n$,
we have:

\begin{theorem}\label{mainthm2}
Let $E$ be a rank $r$ locally-free special sheaf on $Q_n$.
\begin{itemize}
\item If $r\le 2n-1$ and $c_1(E)=0$, then $E$ is semistable;
\item if $r\le n$ and $c_1(E)\ne 0$, then $E$ is stable.
\end{itemize}
\end{theorem}

It is interesting to note that, by \cite[Proposition 4.7]{CMR}, there are no
rank $r\leq n-1$ linear sheaves $E$ on $Q_n$ with $c_1(E)<0$ or
rank $r\leq n-2$ linear sheaves $E$ on $Q_n$ with $c_1(E)=0$.

\begin{proof}
For locally-free special sheaves which are represented as
cohomologies of the monad {\em (M1)}, the statement follows from
Theorem \ref{mainthm} and \ref{mthm2} and for locally-free
special sheaves which are represented as cohomologies of the monad
{\em (M2.1)} and {\em (M2.2)} an analogous argument works.
\end{proof}

Note that using the Fl\o ystad type existence theorem for linear
sheaves on $Q_n$ established in \cite[Proposition 4.7]{CMR}, one
can easily produce examples of rank $2n$ locally-free instanton
sheaves on $Q_n$ as well as rank $n+1$ locally-free linear sheaves
on $Q_n$ which are not semistable, following the ideas in Examples
\ref{ex1} and \ref{ex2}.

However, we do not know whether the bounds in the rank are sharp for
locally-free sheaves on $Q_n$ which are the cohomology of monads of
type {\em (M2.1)} and {\em (M2.2)}. For instance, is there an unstable
rank $2n$ locally-free sheaf on $Q_n$ which can be represented as the
cohomology of a non-linear special monad?


\section{Conclusion}

In this paper we have studied the semistability of torsion-free
sheaves on nonsingular projective varieties with cyclic Picard group
that arise as cohomologies of a particular type of monad. Many
interesting questions regarding linear monads and instanton sheaves
remain unanswered.

First of all, one would like to have a generalizations of Fl\o
ystad's (resp. Costa and Mir\'{o}-Roig's) existence result
\cite{F} (resp. \cite{CMR}) and of Proposition \ref{criterio},
establishing the existence of instanton sheaves over varieties
other than $\pn$ (resp. $Q_n$) and their intrinsic cohomological
characterization.

The semistability of instanton sheaves of low rank indicate the existence
of a well-behaved moduli space of instanton sheaves on cyclic varieties.
One approach to study the moduli space of instanton sheaves would be the
construction of the moduli space of linear monads, using methods from geometric
invariant theory. This task is probably deeply linked with the theory of
representation of quivers, since a monad can be regarded as the representation
of a quiver, the one whose underlying graph is the Dynkyn diagram for $A_3$, into
the category of sheaves, see \cite{K}.

This also brings up the question of a reasonable stability condition for monads,
meaning compatible with geometric invariant theory, and how does it compare with
the slope stability of its cohomology sheaf. Notice that a monad can also be regarded
as an element in the derived category $D^{\rm b}(X)$ of bounded complexes of coherent
sheaves on $X$; the concept of stability on triangulated categories has been recently
introduced by Bridgeland in \cite{Br}, but it is still unclear what does it have to
do with moduli spaces. We hope that the study of the moduli space of instanton sheaves
will shed some light on this topic.


 \end{document}